\input amstex
\input emslab
\magnification=1200
\documentstyle{grgppt01}
\hfuzz=3mm
\parskip=2pt plus 1pt minus 1pt
\footline={\hfil\eightrm \copyright\ Geoffrey Grimmett, 1 July 2002}
\TagsOnRight
\font\smit=cmti8 

\def\es{\varnothing}
\def\ol{\overline}
\def\ii{\itemitem}
\def\E{{\Bbb E}}
\def\P{{\Bbb P}}
\def\cite#1{\citeto{#1}}
\def\section#1{\setcounter{Thm}{0}\bigskip\goodbreak\centerline{{\bf\ignorespaces #1}}
  \nobreak\medskip\nobreak\noindent\ignorespaces}
\def\endsection{\medskip\flushpar}
\def\subsection#1#2{\medskip
  \flushpar{\smc#1\quad#2}\nobreak\smallskip\nobreak\flushpar\ignorespaces
  \nobreak\ignorespaces}

\def\capt#1#2{\botcaption{\tenpoint{\it Table #1.}\q
  #2}\endcaption\par}
\def\mletter#1#2#3{\hskip#2cm\lower#3cm\rlap{$#1$}\hskip-#2cm}
\def\lastletter#1#2#3{\hskip#2cm\lower#3cm\rlap{$#1$}\hskip-#2cm\vskip-#3cm}
\def\figure#1\par{\parindent=0pt
  \vbox{\baselineskip=0pt \lineskip=0pt
  \line{\hfil}
  #1}}

\def\q{\quad}

\def\lf{\lfloor}
\def\rf{\rfloor}
\def\lc{\lceil}
\def\rc{\rceil}

\catcode`\@=11
 \def\logo@{}
\catcode`\@=13
\def\today{\number\day
\space\ifcase\month\or
  January\or February\or March\or April\or May\or June\or
  July\or August\or September\or October\or November\or December\fi
  \space\number\year}

\topmatter
\title Stochastic Apportionment
\endtitle
\author Geoffrey Grimmett
\endauthor
\address
Statistical Laboratory, University of Cambridge, Wilberforce Road,
Cambridge CB3 0WB, United Kingdom
\endaddress
\email g.r.grimmett{\@}statslab.cam.ac.uk\endemail
\http http://www.statslab.cam.ac.uk/$\sim$grg/ \endhttp

\keywords Apportionment problem, stochastic apportionment
\endkeywords
\subjclass 90B99, 91B32, 62C99 
\endsubjclass
\abstract
The problem of how to allocate
to states the seats in the US House of Representatives
is the most studied instance of what is termed the `apportionment problem'.
We propose a new method of apportionment which is stochastic,
which meets the quota condition, and which is fair in the sense of expectations.
Two sources of systematic unfairness are identified, firstly the lower
bound condition (every state shall receive at least one seat),
and secondly the lower quota condition (every state shall
receive at least the integer part of its quota).
\endabstract
\endtopmatter

\document

\section{1. The problem of apportionment}
Ten goats are to be assigned between three brothers in numbers
proportional to the
ages (in years) of the recipients.  
Given the integral nature of a goat, it is not generally
possible to meet exactly the condition of proportionality, and the resulting
`apportionment problem' is a classic of operational research.  The associated
literature is extensive, and includes on the one hand discussions of
criteria to be used in assessing different schemes, and on the other
hand accounts of the properties of specific classes of scheme.  A point
of especial focus has been the apportionment 
of the seats in the House of Representatives
between the states of
the USA.  There are currently
50 states (excluding the District of Columbia)
and 435 seats, which are to be divided between the states
according to the US Constitution [Article I, Section 2]
of 1787 thus: 
$$
\vbox{\hbox{\llap{``}Representatives $\dots$ shall be apportioned 
among the several\hfill}
\hbox{States $\dots$ according to their respective 
Numbers $\dots$'' }}
$$
No scheme is proposed in the Constitution, and the Article 
therein permits a spectrum
of interpretation of the phrase `according to their respective numbers'.
Politicians, lawyers, mathematicians and others have been involved since in the
cyclical debate of how to apportion the seats.

Although there is little in this note which is specific to the US
Congress, we shall, for ease of exposition, use the terminology of the
last problem.  Our targets here are to survey the general area, and to
propose a new method of apportionment which, in a certain way to be made
more precise, meets all the usual criteria for such schemes. This new
method is a lottery scheme whose implementation uses (pseudo-)random
numbers.  The scheme is fair so long as no minimal number of seats is
allocated to all states however small.  In the Congressional example, the
Constitution contains an additional condition that every state shall
receive at least one seat.  It is this violation of the principle of
proportional representation which renders futile all attempts to obtain a
truly fair system in which individuals are equally represented.

We shall in Section 3 describe a new method which meets the so-called quota
condition (see Section 2) and gives proportional representation.
In Section 4 we present an adaptation of this method for use in situations
in which there are lower bounds on the allocations sought. We shall
see in Section 4 that, in the presence of a lower bound, the quota condition
(see Section 2) provides another source of unfairness. 

The established
theory of deterministic apportionment is summarised in Section 5.
Those interested in learning more of the history and practice
of apportionment should consult the excellent book [\cite{BY}]
of Balinski and Young.
\endsection

\section{2. Quota}
We suppose there are $s$ states with respective populations
$\pi_1, \pi_2,\dots, \pi_s$, and that there are $r$ seats in the House of
Representatives.    The total population size is
$\Pi = \pi_1 + \pi_2 +\dots + \pi_s$, and thus  the exact {\it quota\/} of
seats for state $i$ is $q_i = r\pi_i/\Pi$.  The problem is that $q_i$ are not
generally integers, whereas representatives are (by axiom) indivisible.  We
call $\pi = (\pi_1,\pi_2,\dots, \pi_s)$ and $q = (q_1, q_2,\dots, q_s)$ 
the {\it population vector\/} and {\it quota vector\/}, respectively. We
refer to the pair $(\pi,r)$ as a {\it problem\/}.

An {\it allocation\/} is a vector $\alpha=(\alpha_1, \alpha_2,\dots,
\alpha_s)$ of non-negative integers with sum $r$.  An allocation $\alpha$
is said to `satisfy quota' if $\alpha_i \in \{\lf q_i\rf$, $\lc q_i \rc\}$
for all $i$.  An allocation is
said to `violate quota' if it does not satisfy quota.
Here, $\lf x\rf$ denotes the greatest integer not exceeding $x$,
and $\lc x\rc$ denotes the least integer not less than $x$.
We speak of $\lf q_i\rf$ (respectively, $\lc q_i\rc$) as the lower
(respectively, upper) quota for state $i$.

It seems generally (if not universally) accepted that the property of
satisfying quota is desirable.  It is clear that, for all problems of the
above type, there exists necessarily at least one allocation which
satisfies quota.  This can cease to be the case when further conditions are
added.  In a variety of situations including that of the US Congress,
there is a requirement that the $\alpha_i$ be not too small.  Let $l =
(l_1, l_2,\dots, l_s)$ be a given vector of non-negative integers.  We
say that an allocation $\alpha$ `has lower bound' $l$ if $\alpha_i \geq
l_i$ for all $i$.  Of special interest is the case $l=1$, the vector
of ones, not least since the relevant article of the US Constitution 
requires that
$$
\hbox{
``each State shall have at Least one Representative $\dots$''}
$$
Such a requirement is potentially disturbing since there exist
problems $(\pi,r)$ for which no allocation exists
satisfying quota and with lower bound 1.  For a simple (if extreme) example, 
consider the case when $\pi = (1,1,7)$ and $r=3$.

It is a simple matter to see that there exists an allocation which satisfies 
quota and with lower bound $l$ if and only if
$$
l_i \leq \lc q_i\rc  
\text{ for all $i$, and } \sum_i \max\{l_i,
\lf q_i\rf\} \leq r.
$$

The satisfying of quota seems to be regarded as paramount amongst
properties of {\it allocations\/}. As noted above,
we shall see in Section 4
that, in the presence of a lower bound, the requirement of quota can 
be a further source of unfairness.
There are certain desirable properties
of schemes for finding allocations, most prominently that the scheme
avoids each of the three `paradoxes' --- Alabama, population, and
new-state --- to which we shall return in Section 5.

A word for the novice --- in common with other similar problems, one may be
tempted to identify many desirable properties of schemes, only to find that
no scheme has them all.  A well known example of this phenomenon is `Arrow's
Impossibility Theorem', which states that no preference ranking exists for a 
society which embraces a certain collection of five reasonable axioms
(see [\cite{LR}, Chapter 14]).

\endsection

\section{3. Stochastic apportionment without lower bounds}
Although lottery schemes have been mentioned briefly in the literature
(see, for example, [\cite{BY}]), the established theory is
concentrated on
deterministic schemes.  Our purpose here is to propose a family of stochastic
schemes, one in particular, which satisfy quota and
which have the advantage of being truly fair
and proportional in that each state receives a (possibly random) number
of seats having mean value equal to the quota of the state.  We will see
in Section 4 that this cannot generally be achieved in the presence of
lower bounds on allocations, and therefore {\it we make the facilitating
assumption for the duration of this section that no lower bound is
required\/}.

A {\it random allocation\/} is a vector $A=(A_1, A_2,\dots, A_s)$ of 
non-negative-integer-valued random variables with sum $r$. A {\it randomized
scheme} is a mapping which, to each problem $(\pi, r)$, allocates a probability
distribution on the space of appropriate allocations.  Otherwise expressed,
a randomized scheme results in a random allocation (we shall not spend
any time on the choice of probability space, and such like).

There is a subtlety to the notion of a randomized scheme which we discuss
briefly.  We  may seek to apply such a scheme to two 
given problems, perhaps by
applying it twice to the same problem $(\pi, r)$,
or perhaps by applying it to
$(\pi, r)$ and to another problem $(\pi', r')$  obtained from
$(\pi,r)$ by changing some
of the parameters.  In so doing, we encounter the question of `coupling'.
That is, since randomized schemes make use of pseudo-random numbers, we
shall need to specify whether, at the one extreme, we re-use for the
second problem the pseudo-random numbers used already for the first, or, at the
other extreme, we make use of `new' pseudo-random numbers. We discuss this
no further at this stage (see, however, the two
final paragraphs of this section), since the discussion will concentrate for the
moment on the use of randomized schemes for a single problem only.

We say that a random allocation $A$ `satisfies quota almost surely' if 
$$
\P(A
\text{ satisfies quota}) =1. 
$$
Here, $\P$ denotes probability, and (later)  
$\E$ denotes expectation. A randomized scheme is said to satisfy
quota if the ensuing random allocation satisfies quota almost surely.

In advance of making a concrete proposal for a randomized scheme, we point
out that there exist in general a multiplicity of such schemes which satisfy
quota.  Let $q_i' = q_i - \lf q_i \rf$, let $s' =|\{i:q_i' > 0 \}|$, and
let $r' = r - \sum_i \lf q_i\rf$.  We call $q'=(q_i': 1\le i\le s)$
the {\it fractional-quota vector\/}.
Since state $i$ receives $\lf q_i\rf$ seats
of right, we are required only to allocate the remaining $r'$ seats amongst
the $s'$ `unsatisfied' states.  This may be done in $\binom {s'}  
{r'}$ ways, to each of which we must assign a probability.  If we
require in addition that the scheme be {\it fair\/} in the sense that 
the ensuing allocation $A$ satisfies $\E (A_i)=q_i$ for all $i$, then
we obtain thus a set of $s'$ constraints.  Thus the space of such randomized
schemes may have up to $\binom {s'}{r'} -s'$ degrees of freedom.
In the special  cases $r'=1$, $r'=s'-1$, there is a unique randomized scheme
which is fair and satisfies quota.

We turn now to our concrete proposal.  There are three steps.

\ii{I.} We permute at random the labels of the states.  
\flushpar
This is proposed since
the scheme which follows depends on the labelling of the states (that
is, on the indices 
of the $\pi_i$), and it seems desirable to reduce to a minimum any
correlations which depend on this extraneous element.  In the
Congressional example, it would arguably not be right for the allocation
to Alabama to depend systematically 
to a greater degree on that to Alaska than that to
Maine.  {\it For ease of presentation, we shall in the following not
change the notation $\pi_i$, but shall assume that $\pi$ is the population
vector after permuting at random\/}.

\ii{II.}  We provisionally allocate $\lf q_i\rf$ seats to state $i$.

\flushpar
This is prompted by the minimal quota for each state, and leaves so far
unallocated a certain number $r'=r-\sum_i \lf q_i\rf$ of seats.  Let
$q_i' = q_i - \lf q_i\rf$, as before.

\ii{III.}  Let $U$ be a random 
variable with the uniform distribution on $[0,1]$,
and let $Q_i = U+\sum^i_{j=1} q_i'$.  Let $A_i'$ be the indicator function of
the event that the interval $[Q_{i-1}, Q_i)$ contains an integer.  We 
allocate a further $A_i'$ seats to state $i$.

State $i$ receives a total of $\lf q_i\rf + A_i' = A_i$ seats.  The total
number of seats allocated at step III equals the length of the interval
$[Q_0, Q_s)$, which equals $\sum_i \{q_i - \lf q_i\rf\} =r'$, whence $A$ is an
allocation (that is, it has sum $r$).
It is evident that $\lf q_i\rf\le A_i\le \lc q_i\rc$,
whence $A$ satisfies quota. 
In addition, each $Q_j$, when reduced modulo $1$, is uniformly
distributed on the interval $[0,1]$. Therefore,
$\E(A_i')$ equals the length of the interval
$[Q_{i-1},Q_i)$, which is to say that $\E(A_i')=Q_i-Q_{i-1}=q_i'$ and hence
$\E(A_i)=q_i$. We summarise
this by saying that the scheme satisfies quota, and is fair
in the sense that the mean number of seats per head of population is constant
between states.

We offer no justification for this scheme apart from fairness and ease of
implementation.  We caution against adopting any randomized scheme without a
{\it proof\/} of fairness, and as an example we summarise one arguably
reasonable but definitely unfair scheme for allocating the remaining $r'$ seats. 

\noindent
{\bf Conditional sampling.}
We select independent indices $I_t$, $1 \leq t \leq r'$, each $I_t$ having
distribution ${\Bbb P}(I_t=i)=q_i'/\sum_i q_i'$, and 
we consider a random vector
$J=(J_1, J_2,\dots, J_{r'})$ having the probability
distribution of $I_1, I_2, \dots, 
I_{r'}$ conditional on $I_u \neq I_v$ for all $u\neq v$.  We allocate one
extra seat to those states having indices $J_1, J_2,\dots, J_{r'}$.  It is
the case that
$$
\E |\{ t: I_t = i \}| = q_i' \quad\text{for all } i,
$$
but the same statement is generally false with $I$ replaced by $J$.  This 
is easily seen when $r'$ and $s'$ are large.  The conditioning amounts to a 
large deviation, and the distribution of $J$ may be close to 
an appropriately
tilted distribution (see [\cite{GS}, Section 5.11]).

We close this section 
with some remarks on the use of pseudo-random numbers.  It will be argued
that the fairness of a stochastic 
scheme may not be evident to the population. It
seems that politicians may be less willing than the people to accept such
a scheme, since politicians are very sensitive to the marginal value to a
party of a single seat.
Lotteries are already in wide use in areas having impact on individuals,
not least in state-accredited systems for raising money for 
so-called `good causes'.
As a further example, the allocation of individuals to the control group of a
medical trial is usually done by lottery, and such a decision may be a
matter of life or death.  See [\cite{SC}] for an extended discussion of the
drawing of lots.

When applying a randomized scheme to two or more problems,
one needs to decide whether or not to re-sample the required
pseudo-random numbers. The expectations under study remain the same, but the 
external perception of fairness is likely to be greater if the 
roulette wheel is spun afresh. A statistical virtue of re-sampling 
is the reduction
of variances.
\endsection

\section{4.  Stochastic apportionment with lower bounds} 
We consider in
this section the problem of finding an allocation which satisfies quota
and is subject to a lower bound $l$. Let $l=(l_1,l_2,\dots,l_s)$ be a vector
of non-negative integers satisfying $l_i\le \lc q_i\rc$ for all $i$.
We assume in addition that $\sum_i l_i \le r$, since otherwise there
exists no allocation with lower bound $l$ and sum $r$.
We call a state {\it small\/} 
if $q_i<l_i$, and we note
that the lower-bound condition is 
tantamount to favouring the small states.  
It follows that, whenever there exist small states, then the other
states are at a disadvantage, and there can exist no scheme
which is fair.  But how fair a randomized scheme can be achieved?

Let $I_- =\{i: q_i<l_i\}$ be the set of small states, and also
$I_= =\{i:q_i=l_i\}$ and $I_+=\{i:q_i>l_i\}$. We assume that there
exists at least one small state,
in that $I_-\ne\es$, since otherwise there is no new problem.
Since $\sum_i l_i\le r$, it must be the case that $I_+ \ne \es$.
For $i\in I_- \cup I_=$, we allocate to state $i$ exactly
$l_i$ seats, and
we write $\mu=r-\sum_{i\in I_- \cup I_=} l_i$ for the 
number of remaining seats.  Note
that $\sum_{i\in I_+} q_i >\mu$, whence no fair allocation can exist.

In allocating seats to states in $I_+$, we seek a guiding principle, 
and we propose  the principle of
{\it equality of representation amongst 
states in $I_+$\/}.  

The scheme of Section 3 (particularly step III thereof)
is based on the quota vector $q$ and the fractional-quota vector $q'$.
Since the allocations to states in $I_- \cup I_=$ have already been determined,
we are directed by the above principle to seek a new
quota vector $Q=(Q_i: i\in I_+)$ such that:
\comment
fractional-quota vector $q''=(q_i'' : 1\le i\le s)$ such that
$q_i''=0$ for $i\in I_- \cup I_=$, and:
\endcomment
\roster
\item"(a)" $\sum_{i\in I_+} Q_i = \mu$,
\item"(b)" the mean number of representatives per head of population is constant
between states in $I_+$, which is to say that there exists a constant $\gamma$
such that
$$
\frac{Q_i}{q_i} = \gamma \quad \text{for all }
i \in I_+,
\tag 4.1
$$
\item"(c)" $\lfloor q_i \rfloor\le Q_i \le \lceil q_i\rceil$ 
for all $i \in I_+$.
\endroster
Following consideration of the lower bound condition, we shall require also that
\roster
\item"(d)" $Q_i \ge l_i$ for all $i\in I_+$.
\endroster
Note that (c) implies (d) since, by (c),
$$
Q_i\ge \lf q_i\rf \ge l_i\q\text{for } i \in I_+.
$$

We next investigate conditions under which (a), (b), (c) may be achieved
simultaneously.  Assume that (a) and (b) hold.  We sum (4.1) over $i \in I_+$
to obtain by (a) that
$$
\gamma =
\frac{\sum_{i\in I_+} Q_i}{\sum_{i\in I_+} q_i} = 
\frac{\mu}{\sum_{i\in I_+} q_i}.
\tag 4.2
$$
It follows that (c) is satisfied if and only if
$$
Q_i=\gamma q_i \ge \lfloor q_i\rfloor\quad\text{for all } i \in I_+.
\tag 4.3
$$

If (c) holds, which is equivalent to
$\lfloor Q_i\rfloor = \lfloor q_i\rfloor$
for $i\in I_+$,
then we proceed by applying the algorithm of the last section thus.
For $i\in I_+$, we allocate $\lf q_i\rf$ seats to state $i$. Then
we apply step III to the new
fractional-quota vector $Q'=(Q_i-\lfloor Q_i\rfloor: i \in I_+)$ with (by (a))
$\mu - \sum_{i\in I_+}\lf q_i\rf$ seats.  The outcome
is a random allocation which satisfies quota and which is fair when
restricted to the states in $I_+$.

Suppose that (c) does not hold, in that 
there exists $i\in I_+$ such that $Q_i<\lfloor q_i\rfloor$.
One strategy would be to apply the scheme of the last 
section to the quota vector
$Q$, obtaining thereby a random allocation $A$
satisfying $\E(A_i)=Q_i$ for each $i\in I_+$. If this
random allocation does not satisfy quota, then the scheme has failed.
There is however a certain probability that quota is satisfied. By counting
the mean number of states whose allocation 
does not satisfy the respective quota,
we find that
$$
\P(A \text{ does not satisfy quota}) \le 
\sum_{i\in I_+}  \max\{1,\lf q_i\rf-Q_i\} 1_{\{Q_i<\lfloor q_i\rfloor\}},
\tag 4.4
$$
where $1_A$ denotes the indicator function of the event $A$.
Equality holds in (4.4) if there exists at most one state $i$ with $Q_i<
\lfloor q_i\rfloor$.

Let us suppose that (c) does not hold but that
the outcome $A$ satisfies quota. A slightly subtle point is that,
conditional on this event, $A$ is not fair. This is so because the conditioning
changes the expectations. In summary, the scheme --- {\it work with the
quota vector $Q$ repeatedly
until one obtains a random allocation which satisfies quota\/} ---
is not a fair scheme.

A feasible line of enquiry which we have not pursued is to
postulate probabilistic models for problems, and to calculate
the probability for such a model that the above scheme results in
an allocation which does not satisfy quota.  In certain circumstances, the `law
of anomalous numbers' (see [\cite{F2}, \cite{GS}]) could be used
as a basis for such a model.

If the satisfying of quota is paramount, then one must accept
a further degree of unfairness, and we propose the following.
If $Q_i<\lfloor q_i\rfloor$, 
we allocate to state $i$ exactly $\lf q_i\rf$ seats.
Writing $J_+=\{i\in I_+: Q_i\ge\lfloor q_i\rfloor\}$,
we now iterate the above process restricted to the states in $J_+$.
That is, we seek, as above, numbers $R_i$, $i\in J_+$,  such that:
\roster
\item"(a$'$)" $\sum_{i\in J_+} R_i = 
\mu - \sum_{i\in I_+\setminus J_+} \lf q_i \rf$,
\item"(b$'$)" the mean number of representatives 
per head of population is constant
between states in $J_+$, in that there exists a constant $\gamma'$
such that
$$
\frac{R_i}{q_i} = \gamma' \quad \text{for all }
i \in J_+,
$$
\item"(c$'$)" $\lf q_i\rf \le R_i \le \lc q_i\rc$ for all $i \in J_+$.
\endroster
\flushpar
We find this time from (a$'$) and (b$'$) that $R_i=\gamma' q_i$ where
$$
\gamma' =\frac{\mu-\sum_{i\in I_+\setminus J_+} \lf q_i\rf}
 {\sum_{i \in J_+} q_i}.
$$
If $R_i \ge \lf q_i\rf$ for all $i\in J_+$, we apply the scheme to the new 
quota vector $R=(R_i:i\in J_+)$. Otherwise, we note
that any state $i$ with $R_i<\lf q_i\rf$ must be handled
in a way which will disfavour those remaining, namely by allocating to it
$\lf q_i\rf$ seats.

This scheme, when iterated to reach a conclusion, identifies different
levels of unfairness in its consecutive applications of the principle
of {\it equality of representation amongst the remaining states\/}.  Note that it
terminates with a random allocation which satisfies quota if and only if there exists
an allocation which satisfies both quota and the lower bound.

We emphasize that, in the presence of a lower bound, the principle
of satisfying quota is another potential source of
unfairness. The upper quota presents no problem, but
the lower quota can indeed be problematic.

Let us apply the above argument in the Congressional example, with
$l=1$.  In the following table are listed those states $i$ which, in the
nine ten-yearly apportionments of 1920--2000, have $Q_i<\lfloor q_i\rfloor$,
these being the apportionments with the current House size of 435.
In just four of these apportionments, namely
those of 1920/1950/1970/2000, do there exist states $i$
with $Q_i<\lf q_i\rf$, and in each such case
there is only
a small probability that the proposed
scheme results in an allocation based on  
$Q$
which does not satisfy quota.
For example, in 1950, the ensuing allocation could only fail to satisfy
quota if New York were allocated 42 seats rather than 43, 
an event of probability
$\lf q_i\rf-Q_i=43-\text{42.962}=\text{0.038}$, 
where $i$ is
the index of New York.
The corresponding probability
for the year 2000 is $\lf q_j\rf-Q_j=19-\text{18.999}=\text{0.001}$, 
where $j$ is
the index of Pennsylvania. A similar calculation is valid for
the other years, and in each case the resulting probability is small. 
Based on this empirical evidence, we claim that the proposed
scheme, based on $Q$,
is likely to result in an apportionment of the House of Representatives
which satisfies quota.

\vbox{
\medskip\goodbreak
{\tabskip=1em plus2em minus.5em
\halign to \hsize{\tabskip=1em\hfil#\hfil&\hfil#\hfil&\hfil#\hfil&\hfil#\hfil&\hfil#\hfil\tabskip=1em plus2em minus.5em\cr
{\smit Year}&{\eightpoint\it $\#$ small states}&{\eightpoint\it states with $Q_i<\lfloor q_i\rfloor$}&{\eightpoint\it original quota $q_i$}&{\eightpoint\it new quota $Q_i$}\cr
\noalign{\nobreak\vskip3pt\nobreak}
1920&3&Pennsylvania&36.053&35.973\cr
&&Florida&4.004&3.995\cr
1930&3&none&&\cr
1940&3&none&&\cr
1950&3&New York&43.038&42.962\cr
1960&4&none&&\cr
1970&3&Virginia&10.000&9.984\cr
1980&3&none&&\cr
1990&3&none&&\cr
2000&4&Pennsylvania&19.013&18.999\cr
}}

\capt{1}{For the apportionments of 1920--2000, this is
a list of states $i$ for which $Q_i<\lfloor q_i\rfloor$. 
The entries are derived from
data to be found in [\cite{BY}].}
\medskip
}

However, in these years 1920/1950/1970/2000, such an outome would be 
slightly unfair to the other states. Suppose,
instead, we acknowledge this unfairness explicitly by 
automatically allocating to any state $i$ with $Q_i<\lfloor q_i\rfloor$
the number $\lf q_i\rf$ seats. In each case, when we 
calculate the amended quota vector $R=(R_i:i\in J_+)$,
we find that every $R_i$ satisfies $R_i\ge\lfloor q_i\rfloor$. 
The outcome of the composite scheme is therefore an allocation 
which favours firstly
the small states, and secondly any state $i$ with $Q_i<\lfloor q_i\rfloor$, 
and then proceeds to allocate seats to
the remaining states in a way which gives equality of representation
between them.

Finally we note the existence of an alternative method which might 
in principle be followed. This
is easier to apply but violates the condition of
fairness, and is as follows. One reduces the original fractional-quota
vector $(q_i': i\in I_+)$,
by a constant factor, that is, one
works with a new fractional-quota
vector $(q_i'' : i \in I_+)$ given by $q_i'' = \delta q_i'$
satisfying
$$
\sum_{i\in I_+} \{ \lf q_i \rf + q_i''\} = \mu .
$$
\endsection

\section{5. Deterministic apportionment}
It is conventional to apportion the House of Representatives using
an algorithm which is deterministic rather than stochastic.  The book
[\cite{BY}] is an excellent account of theory and practice for
mathematicians and non-mathematicians alike, and the reader is referred to
it for a full account.  
No attempt is made here to do more than to summarise the situation.

Let us assume that there is no lower bound.  The discussion 
in the literature has
concentrated largely on a category of schemes termed `divisor methods'.  Let
$\delta : \{0,1,2,\dots\}\to {\Bbb R}$ be a given function satisfying $b
\leq \delta (b) \leq b+1$ for all $b$.  Let $\lambda$ be a parameter
taking positive values, to be thought of as the notional number of head of
population to be represented by each representative.  The population
$\pi_i$ of state $i$ is divided by $\lambda$ to obtain the 
$\lambda$-{\it quota\/}
$\pi_i/\lambda$ of seats for that state.  At the first stage, we find the  
$\lambda$-{\it allocation\/} to
state $i$, defined as 
$$ 
\alpha_i(\lambda) = \cases \lf \pi_i/\lambda\rf
&\text{if $\pi_i/\lambda \leq \delta (\lf \pi_i/\lambda\rf)$},\\ \lf
\pi_i/\lambda\rf +1 &\text{if $\pi_i/\lambda > \delta (\lf
\pi_i/\lambda\rf)$}. 
\endcases 
$$ 
When equality holds in that $\pi_i/\lambda = \delta (\lf
\pi_i/\lambda \rf)$, we have set $\alpha_i(\lambda) = \lf \pi_i/\lambda
\rf$ for the sake of definitiveness, but it is important only that 
some definite rule be followed. The resulting $\lambda$-allocation
$\alpha(\lambda) = (\alpha_i (\lambda) : 1 \leq i \leq s)$ does not
generally sum to the house size $r$.
At the second stage, we `tune' $\lambda$ until we find $\ol{\lambda}$
such that $\sum_i \alpha_i(\ol\lambda) = r$, and we output the allocation
$\alpha (\ol{\lambda})$.

There is an infinity of possible choices for the function $\delta$, of
which the following instances have been studied.

\medskip
{\tabskip=2em plus2em minus.5em
\halign to \hsize{\tabskip=1.5em\hfil#\hfil&\hfil#\hfil&\hfil#\hfil\tabskip=2em plus2em minus.5em\cr
{\it Name}&{\it Originator}&{\it The function $\delta(b)$}\cr
\noalign{\vskip3pt}
Smallest divisors&Adams (1832)&$b$\cr
Harmonic means&Dean (1832)&$2/\{b^{-1} + (b+1)^{-1}\}$\cr
Equal proportions&Hill (1911)&$\sqrt{b(b+1)}$\cr
Major fractions&Webster (1832)&$b+\tfrac{1}{2}$\cr
Greatest divisors&Jefferson (1792)&$b+1$\cr
}}

\capt{2}{The five principal divisor methods ordered by increasing $\delta$.}
\medskip

A sixth scheme, termed the method of Hamilton (1792)
or the `method of largest
remainders', is as follows. A state with quota $q_i$ receives by right
$\lf q_i\rf$ seats.  The remaining $r-\sum_i \lf q_i\rf$ seats are
allocated to those $r-\sum_i \lf q_i\rf$ states $j$ 
with largest remainders $q_j - \lf q_j\rf$.

One seeks principles which enable distinctions to be drawn between these
six schemes.  In common with other instances in operational
research, one can be over-principled. Every method has its
drawbacks, and to seek the perfect system can be to eliminate all
possibilities.  A detailed and informative discussion is to be found in
[\cite{BY}], from which a few points are extracted here.

A scheme is called {\it population monotone\/} if: for all $i\neq j$, when
individuals move from state $i$ to state $j$, then $i$ should not get more
seats {\it and\/} $j$ fewer.  We have from [\cite{BY}, Thm 6.1] that no
scheme exists which is both population monotone and necessarily satisfies
quota.  Furthermore, (see [\cite{BY}, Prop.\ 6.4]), Jefferson's `greatest
divisor' method is the only population monotone method which stays above
lower quota (in that the ensuing allocation $\alpha$ necessarily satisfies
$\alpha_i \geq \lf q_i\rf$ for all $i$), and Adams' `smallest divisor'
method is the only population monotone scheme which stays below upper
quota (in that $\alpha_i \leq \lc q_i\rc$ for all $i$).

A point of focus is the degree to which a scheme favours large over small
states, There are various ways of measuring such bias, both theoretical and
empirical, and the case is made in [\cite{BY}, \cite{Y}] that
Webster's `major fraction' method is the least biased in this regard.

A scheme is said to be `house monotone' if no state's allocation
diminishes when the size $r$ of the house increases.  A scheme which is
not house monotone is said to suffer from the `Alabama paradox'. It is
considered desirable that a scheme be 
house monotone.  Other desirable features of
schemes include the absence of what are known as the
`population paradox' (in a growing population, state $i$ can grow faster
than state $j$, and yet lose a seat to $j$), and the `new-states paradox'
(a new state may join the union, with an appropriate number of new seats,
but the allocations between the original states change).

One may ask in what sense does the stochastic apportionment scheme of Section 3
meet these requirements. The latter scheme is population monotone in a
stochastic sense, which is to say that, when state $i$ loses people to
state $j$, the number of seats allocated to $i$ (respectively, $j$) is
stochastically non-increasing (respectively, non-decreasing). (See
[\cite{GS}, Section 4.12] for a definition of stochastic ordering.) In a
similar stochastic sense, the scheme is house monotone. 

Finally we consider the case of deterministic schemes in situations 
where there is a 
non-trivial lower
bound on the allocation sought. As said already, there may exist no allocation
which satisfies quota {\it and\/} the lower bound, 
and even if there exists such an allocation, 
there will generally exist
no scheme which is fair across the board. The method used currently for
apportioning the House of Representatives is to allocate one seat
to each state however small (in 2000 there were just four states whose quotas
were smaller than one) and then to apply the method of equal proportions to
those states whose `residual' quotas are strictly positive. The outcomes of
this scheme have, fortunately in recent decades, been allocations which
have satisfied quota. See [\cite{BY}, Appendix B] for further information.
\endsection

\bigskip\goodbreak\flushpar
{\bf Acknowledgements.} 
The author thanks Graeme Milton for telling him of the apportionment
problem, and Richard Weber for discussing it with him and for help with
related references. Of great value have been certain documents available on the
web, including writings of Lawrence R.\ Ernst, H.\ Peyton Young, and Roman Shapiro.

\endsection

\enddocument\vfill\eject
\end